\documentclass[reqno]{siamart220329}
\usepackage[english]{babel}
\usepackage{lipsum}
\makeatletter

\newcommand{\authorfootnotes}{\renewcommand\thefootnote{\@fnsymbol\c@footnote}}%
\makeatother
\usepackage{etoolbox} %

\newcommand*\linenomathpatch[1]{%
  \cspreto{#1}{\linenomath}%
  \cspreto{#1*}{\linenomath}%
  \csappto{end#1}{\endlinenomath}%
  \csappto{end#1*}{\endlinenomath}%
}

\linenomathpatch{equation}
\linenomathpatch{gather}
\linenomathpatch{multline}
\linenomathpatch{align}
\linenomathpatch{alignat}
\linenomathpatch{flalign}

\usepackage{amsfonts}
\usepackage{amssymb}
\usepackage{centernot}
\usepackage{graphicx}
\usepackage{algorithm}
\usepackage{algorithmic}
\usepackage[algo2e]{algorithm2e}
\newtheorem{thm}{Theorem}[section]
\newtheorem{prop}[thm]{Proposition}

\usepackage{fullpage}
\usepackage{lineno}

\usepackage{url}

\usepackage{color}

\title{Solving Fokker-Planck equations using the zeros of Fokker-Planck operators and the Feynman-Kac formula}

\begin{document}


\maketitle

 {\normalsize
 \centering
  \authorfootnotes
  Pinak Mandal\footnote{\thanks{Corresponding author: \texttt{pinak.mandal@icts.res.in}}}\textsuperscript{1,3}, Amit Apte\textsuperscript{2,3}\par
  \textsuperscript{1} The University of Sydney, NSW 2006 Australia \par
  \textsuperscript{2}Indian Institute of Science Education and Research, Pune 411008 India\par 
  \textsuperscript{3} International Centre for Theoretical Sciences - TIFR, Bangalore 560089 India \par
  \bigskip}

\begin{abstract}
First we show that physics-informed neural networks are not suitable for a large class of parabolic partial differential equations including the Fokker-Planck equation. Then we devise an algorithm to compute solutions of the Fokker-Planck equation using the zeros of Fokker-Planck operator and the Feynman-Kac formula. The resulting algorithm is mesh-free, highly parallelizable and able to compute solutions pointwise, thus mitigating the curse of dimensionality in a practical sense. We analyze various nuances of this algorithm that are determined by the drift term in the Fokker-Planck equation. We work with problems ranging in dimensions from 2 to 10. We demonstrate that this algorithm requires orders of magnitude fewer trajectories for each point in space when compared to Monte-Carlo. We also prove that under suitable conditions the error that is caused by letting some trajectories (associated with the Feynman-Kac expectation) escape our domain of knowledge is proportional to the fraction of trajectories that escape.
\end{abstract}

\section{Introduction}
\label{sec-intro}
From the motion of a particle suspended in a fluid \cite{karatzas1991brownian},  enzyme kinetics \cite{allen2010introduction} to dynamics of a stock price \cite{karoui1997non}, \cite{delong2013backward} evolving systems in the real worlds are often modelled as systems of ordinary differential equations propagating under the influence of additive noise. These models known as stochastic differential equations (SDE) \cite{oksendal2003stochastic}, \cite{gardiner2009stochastic}, \cite{strauss2017hitch} are directly linked to Fokker-Planck equations (FPE) \cite{risken1996fokker} or Kolmogorov forward equations that describe the evolution of probability density of the state vector. In a prequel \cite{mandal2023learning} we developed a deep learning algorithm to compute non-trivial zeros of high-dimensional Fokker-Planck operators in a mesh-free manner. In this paper we will devise an algorithm to compute solutions of high-dimensional time-dependent FPEs in a mesh-free manner. We will begin by noting an algorithm similar to the one used in \cite{mandal2023learning} to solve stationary FPEs (SFPE) fails for time-dependent FPEs. A widely adopted strategy for solving high-dimensional PDEs is to appeal to Feynman-Kac type formulae \cite{del2004feynman}, \cite{jefferies2013evolution}, since they allow pointwise calculation of solutions without requiring a mesh thus mitigating at least one aspect of the curse of dimensionality. For example, Kakutani's solution of Dirichlet problem for the Laplace operator \cite{kakutani1944131}, \cite{kakutani1944143}, Muller's walk-on-spheres method for Dirichlet problems \cite{muller1956some} and an analogous method called walk-on-stars for Neumann problems \cite{sawhney2023walk}, multi-level Picard iteration method for solving semilinear heat equations \cite{hutzenthaler2021multilevel} are all based on Feynman-Kac type formulae. In recent times  deep learning methods have been combined with the Feynman-Kac formula to solve high-dimensional PDEs \cite{han2018solving}, \cite{blechschmidt2021three}. Even though FPEs are semilinear, parabolic PDEs whose solutions are probability densities and Deep-BSDE method proposed in \cite{han2018solving} deals with semilinear, parabolic PDEs, generic FPEs pose many challenges that make them unapproachable for deep-BSDE. Non-Lipschitzness of drift functions leading to blow-up of SDE trajectories \cite{chow2014almost}, \cite{li2011lack} and unboundedness of the divergence of drift functions causing FPEs to dissatisfy one of the requirements for the Feynman-Kac formula are foremost amongst these challenges. In this paper we apply the Feynman-Kac formula on an auxiliary equation and combine the solution with the zero of the Fokker-Planck operator obtained through the method described in \cite{mandal2023learning} to produce the normalized solution to the time-dependent Fokker-Planck equation.  We will apply our method for problem dimensions ranging from $2$ to $10$ to verify its effectiveness in high dimensions.

As we have noted, real world systems are often modelled as SDEs and we often observe such systems partially due to limited resources and are tasked with determining the distribution of the state vector at a certain time given all the observations up to that time. This is known as the filtering problem in the field of data assimilation and is useful for a variety of topics - global positioning system, target tracking, monitoring infectious diseases, to name a few \cite{sarkka2023bayesian}. FPEs are naturally connected to data assimilation when the underlying dynamics is stochastic. We will show how this method can be used to to calculate the one-step filtering density in the nonlinear filtering problem. To that end we will focus on systems with attractors since such systems are often used as important test cases in the field of data assimilation \cite{carrassi2022data}.

\section{Problem Statement}
\label{sec-prob}
\subsection{Time-dependent FPEs}
Consider the Itô SDE with coefficients $\mu\in C^1(\mathbb R^d;\mathbb R^d)$ and $\sigma\in C(\mathbb R^d; \mathbb R^{d\times l})$ that are independent of time,
\begin{equation}
\begin{aligned}
    dX_t &= \mu(X_t)\,dt + \sigma(X_t)\,dW_t\\
    X_0&\sim p_0\label{eq:SDE-0}
\end{aligned}
\end{equation}
such that $D=\frac{1}{2}\sigma\sigma^\top$ is a positive definite matrix for all inputs. Such an SDE is directly related to the Fokker-Planck equation describing the evolution of the probability density of $X_t$ \cite{risken1996fokker}, \cite{bogachev2022fokker}, \cite{bris2008existence},
\begin{equation}
\begin{aligned}
    &\frac{\partial p}{\partial t} = \mathcal Lp \stackrel{\rm def}{=}
  -\nabla \cdot(\mu p) + {\rm tr}(D\odot\nabla^2p)=0,\qquad\mathbf x\in\mathbb R^d,\; t\in(0, T]\\
  &p(0, \mathbf x) = p_0(\mathbf x),\qquad\mathbf x\in\mathbb R^d\\
  &\int_{\mathbb R^d}p(t, \mathbf x)\,d\mathbf x = 1,\qquad t\in[0, T] 
\label{eq:FPE-0}
\end{aligned}
\end{equation}
where $\odot$ is the Hadamard product and $\nabla^2$ denotes Hessian. The operator $\mathcal L$ is known as the Fokker-Planck operator (FPO). Our goal is to solve \eqref{eq:FPE-0} in a mesh-free way in high-dimensions. In our examples we stick to $\sigma \equiv cI_d$ for some $c>0$ which lets us abuse notation and use $\sigma$ and $D$ as scalar quantities, note that however our final method can be easily applied to problems where $\sigma$ is indeed a function of space or non-diagonal with minor modifications. With this simplification \eqref{eq:FPE-0} becomes,
\begin{align}
    \frac{\partial p}{\partial t}=
  -\nabla \cdot(\mu p) + D\Delta p=0\label{eq:FPE-1}
\end{align}
where $\Delta$ is the Laplacian operator. An extensive amount of work has been done to find numerical solutions to FPEs over the years, for a brief overview the reader can see section 4 of \cite{mandal2023learning}.  
\subsection{One step filtering problem}\label{ssec-1-filter-prob}
Suppose we partially observe $X_t$ at discrete times $t=g, 2g,3g,\cdots$ with observation gap $g$. These observations are given by, 
\begin{align}
    y_k  = Hx_k + \eta_k,\qquad k=1,2,\cdots\label{eq:filtering-obs}
\end{align}
where $x_k = X_{kg}$, $H:\mathbb R^d\to\mathbb R^q$ is a projection matrix and $\eta_k\sim\mathcal N(0_q, \sigma_o^2I_q)$ are iid Gaussian errors in observation. Given all the observations up to time $t=gk$, the filtering problem asks us to compute the distribution of the state vector at that time i.e. $p(x_k|y_{1:k})$. Here, as a simple application, we will calculate the one-step filtering density $p(x_1|y_1)$.

\section{Examples}
\label{sec-example}
All the examples here parallel the examples that appear in the prequel \cite{mandal2023learning} where we call a system a \textit{gradient system} if the corresponding $\mu$ can be written as the gradient of some potential function,
\begin{align}
    \mu=-\nabla V\label{eq:grad-mu}
\end{align}
and otherwise we call it \textit{non-gradient system}. Gradient systems provided important test cases in the prequel \cite{mandal2023learning} while solving stationary FPEs since their steady states are known in analytical form. But for time-dependent FPEs analytical solutions are not known in general even if $\mu$ satisfies \eqref{eq:grad-mu}. 

\subsection{Gradient systems} We use the following gradient systems to verify the effectiveness of our algorithm in high-dimensions. A corresponding non-trivial zero of the Fokker-Planck operator was calculated in \cite{mandal2023learning} for each case.
\subsubsection{2D ring system}\label{ssec-2D}
For $V=(x^2+y^2-1)^2$ and $\mu=-\nabla V$ we get the following FPE,
\begin{align}
     \frac{\partial p}{\partial t}=4(x^2+y^2-1)\left(x\frac{\partial p}{\partial x}+y\frac{\partial p}{\partial y}\right) + 8(2x^2+2y^2-1)p + D\Delta p\label{eq:ring2D}
\end{align}
The corresponding ODE system has the unit circle as a global attractor. We solve this system for $D=1$. As the initial condition we use an equal Gaussian mixture with the components having centers at $\left(-\frac{1}{2},-\frac{1}{2} \right)$ and $\left(\frac{1}{2},\frac{1}{2} \right)$ with covariance matrix $\frac{1}{4}I_2$.
\begin{align}
    p_0(x, y) = \frac{1}{\pi}\exp\left(-2\left(x+\frac{1}{2}\right)^2-2\left(y+\frac{1}{2}\right)^2\right)+\frac{1}{\pi}\exp\left(-2\left(x-\frac{1}{2}\right)^2-2\left(y-\frac{1}{2}\right)^2\right)\label{eq:p0-2D}
\end{align}
This system has a unique solution, for a proof see appendix~\ref{ssec-unique}.
\subsubsection{2nD ring system}\label{ssec-2nD} We can daisy-chain the previous system to build decoupled systems in higher dimensions. In this case the potential is given by
\begin{align}
    V(\mathbf x)=\sum_{i=0}^{\frac{d}{2}-1}(x_{2i}^2+x_{2i+1}^2-1)^2,\qquad d=2n\label{eq:potential-2nD}
\end{align}
We use the following initial condition which can be obtained by daisy-chaining the initial condition in \eqref{eq:p0-2D}.
\begin{align}
    p_0(\mathbf x) = \prod_{i=0}^{\frac{d}{2}-1
    }\left[\frac{1}{\pi}\exp\left(-2\left(x_{2i}+\frac{1}{2}\right)^2-2\left(x_{2i+1}+\frac{1}{2}\right)^2\right)+\frac{1}{\pi}\exp\left(-2\left(x_{2i}-\frac{1}{2}\right)^2-2\left(x_{2i+1}-\frac{1}{2}\right)^2\right)\right]\label{eq:p0-2nD}
\end{align}
Since our algorithm does not differentiate between coupled and decoupled systems, this example serves as a great high-dimensional test case.  Although the analytical solution for this problem is not known, the decoupled nature of this problem implies that we can compare our solution to the 2nD ring system with the solution to the 2D ring system which can be easily computed with other methods (e.g. Monte Carlo) that are not efficient in higher dimensions. Here we solve this system for $n=1,2,3,4,5$ and $D=1$. This system has a unique solution which follows directly from the fact that the 2D ring system has a unique solution.
\subsection{Non-gradient Systems}
Even though the analytic solution for the gradient case is not known for time-dependent FPEs, the gradient structure of the drift makes solving time-dependent FPEs much \textit{easier} as we will see in section~\ref{ssec-limit}. Non-gradient systems on the other hand pose a much harder challenge due to the blow-up of auxiliary SDEs, see section~\ref{ssec-limit} for more details. Here we deal with the following non-gradient systems. A corresponding non-trivial zero of the Fokker-Planck operator was calculated in \cite{mandal2023learning} for each case.
\subsubsection{Noisy Lorenz-63 system}\label{ssec-63} The Lorenz-63 system was first proposed by Edward Lorenz 
as an oversimplified model for atmospheric convection \cite{lorenz1963deterministic}. The corresponding ODE possesses the famous butterfly attractor.
This system and its variants like Lorenz-96 are staple test problems in the field of data assimilation \cite{carrassi2022data}, \cite{yeong2020particle} which is why we also use this system for the calculation of one-step filtering density. We use the standard parameters to define the drift.
\begin{equation}
\begin{aligned}
    &\mu=[\alpha (y-x),\, x(\rho-z) - y,\, xy - \beta z]^\top\\
    &\alpha = 10 \,, \, \beta = \frac{8}{3}\,, \, \rho=28 \
\end{aligned}\label{eq:mu-L63}
\end{equation}
We solve this system for $D=50$ with the following Gaussian mixture as the initial condition,
\begin{align}
    p_0(\mathbf x)=\frac{1}{2\sqrt{(2\pi)^3}}\exp\left(-\frac{\|\mathbf x + 2\times\mathbf 1_3\|_2^2}{2}\right) + \frac{1}{2\sqrt{(2\pi)^3}}\exp\left(-\frac{\|\mathbf x - 2\times\mathbf 1_3\|_2^2}{2}\right)\label{eq:p0-L63}
\end{align}
where $\mathbf 1_3$ denotes the vector with all entries $1$ with respect to the standard basis in $\mathbb R^3$. This system has a unique solution, for a proof see appendix~\ref{ssec-unique}.
\subsubsection{Noisy Thomas system}\label{ssec-Thomas}
The deterministic or ODE version of this system was proposed by René Thomas~\cite{thomas1999deterministic}. This is a 3-dimensional system with cyclical symmetry in $x, y, z$ and the corresponding ODE system has a cyclically symmetric strange attractor.
\begin{equation}
\begin{aligned}
    &\mu=[\sin y - bx,\, \sin z - by,\, \sin x - by ]^\top\\
    &b =  0.2 
\end{aligned}\label{eq:mu-Thomas}
\end{equation}
We solve this system for $D=1$ with the following Gaussian mixture as the initial condition,
\begin{align}
    p_0(\mathbf x)=\frac{1}{2\sqrt{(2\pi)^3}}\exp\left(-\frac{\|\mathbf x + 2\times\mathbf 1_3\|_2^2}{2}\right) + \frac{1}{2\sqrt{(2\pi)^3}}\exp\left(-\frac{\|\mathbf x - 2\times\mathbf 1_3\|_2^2}{2}\right)\label{eq:p0-Thomas}
\end{align}
This system has a unique solution, for a proof see appendix~\ref{ssec-unique}.

To verify our method we will compare the solutions obtained by our method with Monte Carlo simulations of \ref{eq:SDE-0}, for details of the Monte Carlo algorithm see appendix~\ref{ssec-algo-MC}.

\section{The \texorpdfstring{$p_\infty+$}{Lg}FK algorithm}
\label{sec-algo}
First, we will go through the main challenges for solving \eqref{eq:FPE-0} in high dimensions. These challenges and their solutions will naturally lead us to an algorithm for solving high-dimensional FPEs. Lastly, we will discuss various nuances associated with this algorithm.

\subsection{Failure of the physics-informed way}
 Since we used deep learning to find zeros of the Fokker-Planck operator in the prequel \cite{mandal2023learning}, the most natural question becomes does an analogous algorithm work for time-dependent FPEs? For an overview of learning solutions to PDEs in a \textit{physics-informed} fashion see section 5 of \cite{mandal2023learning} or \cite{raissi2019physics}, \cite{blechschmidt2021three}, \cite{sirignano2018dgm}. 
As noted in section 5 of \cite{mandal2023learning} the first step to learning the solution to a PDE is to convert it into an optimization problem. For notational convenience let us first define the function space $\mathcal F$ as follows,
\begin{align}
\mathcal F\stackrel{\rm def}{=}\{&f:[0, T]\times\mathbb R^d\to[0, +\infty):\int_{\mathbb R^d}f(t, \mathbf x)\,d\mathbf x=1\;\forall\; t\in[0, T]\}\label{eq:def-fn-space}
\end{align}

Proceeding parallely to section 5 of \cite{mandal2023learning} or \cite{sirignano2018dgm} we might be tempted to solve the following optimization problem,
\begin{align}
    \underset{f\in\mathcal F }{\rm arg\,inf}\;\mathcal J(f)\stackrel{\rm def}{=}\underset{f\in\mathcal F }{\rm arg\,inf}\left[\frac{1}{MN}\sum_{i=1}^M\sum_{j=1}^N\left(\frac{\partial f}{\partial t}(t_i, \mathbf x_j)-\mathcal Lf(t_i, \mathbf x_j)\right)^2+ \frac{1}{N}\sum_{j=1}^N(f(0, \mathbf x_j)-p_0(\mathbf x_j))^2\right]\label{eq:dynamic-opt}
\end{align}
where $\{\mathbf x_j\}_{j=1}^N$ is a uniform sample from our compact  domain of interest $\Omega_I$ (see section 6.1 in \cite{mandal2023learning}) and $\{t_i\}_{i=1}^M$ is a sample from $[0,T]$.
But it turns out problem \eqref{eq:dynamic-opt} is not a \textit{well-behaved} problem. The following proposition explains why this is the case.
\begin{prop}
    Suppose \eqref{eq:FPE-0} has a unique, strong solution $p$ and 
    \begin{align}
        \mathcal Lf = 0\label{eq:SFPE}
    \end{align}
    also has a unique, strong, probability solution $p_\infty$. 
    Then there exists a sequence of functions $f_k\in\mathcal F$ independent of the samples $\{\mathbf x_j\}_{j=1}^N, \,\{t_i\}_{i=1}^M$ such that \begin{align}
        \lim_{k\to+\infty}\mathcal J(f_k) = 0
\end{align}
    but
    \begin{align}
        \lim_{k\to+\infty}f_k\neq p
    \end{align}\label{prop:failure}
\end{prop}
\begin{proof}
WLOG assume that $\{t_i\}_{i=1}^M$ is ordered and $0=t_1<t_2<\cdots<t_M$. Since we have included $t_1=0$ in our sample, while considering $\frac{\partial}{\partial t}$ at $t=t_1=0$ we'll only consider the right derivative.
For $k\in\mathbb N$ define the sequence,
\begin{align}
    f_{k}(t, \mathbf x) =\begin{cases}
    &p(t,\mathbf x),\;\;t\le\frac{1}{k}\\
    &\phi_k(t)\,p(t, \mathbf x) + (1-\phi_{k}(t))\,p_\infty(\mathbf x),\;\; t>\frac{1}{k}
    \end{cases}
\end{align}
where,
\begin{align}
    \phi_{k}(t) = e^{1-kt}
\end{align}
Since $p,p_\infty$ are probability densities, it is easy to see that  $f_k\in\mathcal F,\;\;\forall\, k\in\mathbb N$.
The second term or the initial condition term in $\mathcal J(f_k)$ vanishes since $f_k(0,\mathbf x)=p(0,\mathbf x)=p_0(\mathbf x)$. Moreover since $p$ is a zero of $\frac{\partial}{\partial t}-\mathcal{L}$, 
\begin{align}
    \sum_{j=1}^N\left(\frac{\partial f_k}{\partial t}(t_1, \mathbf x_j)-\mathcal Lf_k(t_1, \mathbf x_j)\right)^2 =0
\end{align}
where as noted before we have only considered the right time-derivative at  $t=t_1=0$. Now pick $K>0:\frac{1}{K}<t_2$. For $k>K$, using the fact that $p, p_\infty$ are both zeros of $\frac{\partial}{\partial t}-\mathcal{L}$ we see that,
\begin{align}
    \mathcal J(f_{k})=&\frac{1}{MN}\sum_{i=2}^M\sum_{j=1}^N\left(\frac{d \phi_k}{d t}(t_i)(p_\infty(\mathbf x_j)-p(t_i, \mathbf x_j))\right)^2
\end{align}
Due to continuity of $p,p_\infty$ and compactness of $\Omega_I$ there exists $B\ge0:$ 
    \begin{align}
        \sup\left\{p(t, \mathbf x)^2, p_\infty(\mathbf x)^2 :t\in[0, T],\;\mathbf x\in\Omega_I\right\}= B
    \end{align}
which implies $\;\forall\;k>K$,
\begin{align}
    \mathcal{J}(f_k)\le&\frac{2B}{MN}\sum_{i=2}^M\sum_{j=1}^N\left(\frac{d \phi_k}{d t}(t_i)\right)^2=\frac{2Be^2k^2}{MN}\sum_{i=2}^M\sum_{j=1}^Ne^{-2kt_i}
<{2Be^2k^2}e^{-2kt_2}
\end{align}
Since $t_2>0$ we have
\begin{align}
    \lim_{k\to+\infty}\mathcal J(f_k)=0
\end{align}
But clearly the pointwise limit of $f_k$ is given by
\begin{align}
    f(t, \mathbf x) =\begin{cases}
    &p_0(\mathbf x),\;\;t=0\\
    &p_\infty(\mathbf x),\;\; t>0
    \end{cases}\label{eq:limit-min-seq}
\end{align}
\end{proof}
The proof assumes very little about the operator $\mathcal L$ and hence can be extended to other partial differential equations of parabolic type. Only assumptions required to prove proposition~\ref{prop:failure} are the existence of strong, unique probability solutions which are true for all of our example problems, see appendix 9.1 in \cite{mandal2023learning} and appendix~\ref{ssec-unique} in this paper. Since we only have finite computational resources it is sensible to investigate loss functions with finite sample-sizes and note that the proof of proposition~\ref{prop:failure} works for any choice of samples $\{t_i\}_{i=1}^M,\,\{\mathbf x_j\}_{j=1}^N$ with arbitrary but finite sample-sizes. Note that the pathological sequence we constructed is independent of the samples $\{t_i\}_{i=1}^M,\,\{\mathbf x_j\}_{j=1}^N$ which is an important component of the failure of the physics-informed way. In the case of stationary Fokker-Planck equations we can also construct pathological functions that make the physics-informed loss (defined in section~6.3 of \cite{mandal2023learning}) vanish without converging to a solution by simply interpolating a solution at the sample points with non-solutions at unsampled points. But the physics-informed method still produces correct solutions in the stationary case since these pathological functions are sample-dependent and as the sample-size increases, they resemble a solution more and more closely. From the perspective of training a neural network to optimize problem~\ref{eq:dynamic-opt} the main difficulty is revealed by the limit of the minimizing sequence given in \eqref{eq:limit-min-seq}. The network can behave like the initial condition for a short amount of time and then mimic $p_\infty$ for all subsequent times losing all variation in time. Moreover, one can construct many such sequences with pathological behaviors, not to mention restricting neural networks to $\mathcal{F}$ and have them be expressive is another challenge since normalization is an extremely difficult condition to implement in high dimensions, see the discussion in sections 2 and 6.4 of \cite{mandal2023learning}. Apart from these difficulties, for other explorations of failure modes of physics-informed neural networks see \cite{krishnapriyan2021characterizing}, \cite{basir2022investigating}.
\label{ssec-failure}

\subsection{The Feynman-Kac formula}
Due to the ineffectiveness of the physics informed method, we turn to another very useful paradigm for solving high dimensional PDEs. The Feynman-Kac formula \cite{pham2015feynman} has been used very successfully in recent years to solve high-dimensional PDEs \cite{hutzenthaler2021multilevel}. They have also been used alongside deep-learning to create algorithms for semilinear parabolic PDEs \cite{han2018solving}. Moreover, while verifying algorithms for very high dimensional PDEs one often relies on the Feynman-Kac formula to get a reference solution, see for example the examples given in sections 4, 5 of \cite{sirignano2018dgm}. The pointwise nature of the Feynman-Kac formula makes it very attractive for high-dimensional problems, since it waives the requirement for a mesh. Below we briefly discuss a few versions of the formula and the challenges associated with using it for \eqref{eq:FPE-0}.

The following version of the formula appears in \cite{oksendal2003stochastic} and it is one of the most well-known versions.
\begin{thm}[Feynman-Kac]
Suppose $f\in C^2_0(\mathbb R^d)$ and $q\in C(\mathbb R^d)$. If $q$ is bounded below then
\begin{align}
    v(t, \mathbf x) = \mathbb E\left[\left.\exp\left(-\int_0^t q(\bar X_s)\,ds\right)f(\bar X_t)\right\vert \bar X_0=\mathbf x\right]\label{eq:FK-sol}
\end{align}
is the unique solution of
\begin{equation}
\begin{aligned}
    &\frac{\partial u}{\partial t}= Au - qu\\
    & u(0, \mathbf x)=f(\mathbf x)
\end{aligned}\label{eq:FK}
\end{equation}
that's bounded on $K\times\mathbb R^d$ for every compact $K\subset[0,+\infty)$, where $\bar X_t$ is a time-homogeneous Itô process satisfying 
\begin{align}
    d\bar X_t = \bar\mu(\bar X_t)\, dt + \bar\sigma(\bar X_t)\,dW_t\label{eq:ito-time-homo}
\end{align}
with
\begin{align}
    |\bar\mu(\mathbf x)-\bar\mu(\mathbf y)| + |\bar\sigma(\mathbf x)-\bar\sigma(\mathbf y)|\le C|\mathbf x -\mathbf y|,\;\forall\,\mathbf x, \mathbf y\in\mathbb R^d\label{eq:ito-lipschitz}
\end{align}

where $C$ is a positive constant and $A$ is the infinitesimal generator of $\bar X_t$ given by 
\begin{align}
    Af = \bar\mu\cdot\nabla f + \frac{1}{2}{\rm tr}(\bar\sigma\bar\sigma^\top\odot\nabla^2f)\label{eq:ito-gen}
\end{align}
\label{thm:FK-Oksendal}\end{thm}

Note the three requirements for theorem~\ref{thm:FK-Oksendal} are global Lipschitzness of $\bar\mu$ and $\bar\sigma$, vanishing of the initial condition $f$ at infinity and the lower-boundedness of $q$. Lipschitzness of drift and diffusion coefficient ensures that SDE~\eqref{eq:ito-time-homo} has a unique, non-exploding, strong solution \cite{ikeda2014stochastic}. Vanishing of the initial condition $f$ and lower-boundedness of $q$ guarantee that the expectation in \eqref{eq:FK-sol} is bounded.

If we try to apply formula \eqref{eq:FK-sol} directly to the Fokker-Planck equation \eqref{eq:FPE-1}, we have to identify the following quantities,

\begin{equation}
\begin{aligned}
    &\bar\mu = -\mu\\
    &\bar\sigma =\sigma\\
    & q = -\nabla\cdot\mu
\end{aligned}
\end{equation}

In our examples none of the three requirements for theorem~\ref{thm:FK-Oksendal} are always satisfied. Concretely, the 2D ring system \eqref{eq:ring2D} satisfies none of the requirements. Lipschitzness of the drift is too restrictive a condition to model real world phenomena and therefore attempts have been made to relax this condition. In \cite{yamada1971uniqueness} the Lipschitz condition is replaced with Hölder condition for 1D SDEs, furthermore in \cite{fu2010stochastic}, \cite{li2011strong} sufficient conditions for existence of pathwise unique solutions of 1D SDEs have been discussed. \cite{xi2019jump} discusses the more general case when the SDEs are multi-dimensional with non-Lipschitz drifts and gives us an identical Feynman-Kac formula for \eqref{eq:FK} under the condition that \eqref{eq:ito-time-homo} has a unique, weak solution. For the sake of completeness, below we provide a special case of the more general result proved in theorem~7.5 of \cite{xi2019jump}.

\begin{thm}[Feyman-Kac with non-Lipschitz drift] Assume the following.
\begin{enumerate}
    \item \eqref{eq:ito-time-homo} has a unique weak solution.
    \item $v(\cdot,\cdot):[0, T]\times\mathbb R^d\to\mathbb R$ lies in $C^{1,2}([0, T)\times\mathbb R^d)\cap C_b([0, T]\times\mathbb R^d)$ and satisfies the Cauchy problem \eqref{eq:FK}.
    \item $f$ is uniformly bounded.
\end{enumerate}
Then $v$ is given by the formula \eqref{eq:FK-sol}.
\end{thm}

The requirement that $q$ be lower-bounded can be easily arranged while simultaneously simplifying \eqref{eq:FK-sol} by looking at an auxiliary equation. To derive this auxiliary equation we can write, 
\begin{align}
    p(t, \mathbf x) = h(t, \mathbf x)p_\infty(\mathbf x)
\end{align}
where $p_\infty$ is a nonzero zero of the operator $\mathcal L$. It is then easy to see that $h$ satisfies the equation, 
\begin{equation}
\begin{aligned}
    &\frac{\partial h(t, \mathbf x)}{\partial t} = (\sigma^2\nabla\log p_\infty - \mu)\cdot \nabla h + \frac{\sigma^2}{2}\Delta h,\qquad \mathbf x\in \mathbb R^d, t\in(0, T]\\
    &h(0, \mathbf x) = \frac{p_0(\mathbf x)}{p_\infty(\mathbf x)},\qquad \mathbf x\in \mathbb R^d
\end{aligned}\label{eq:h}
\end{equation}
Note that $\mathcal L$, being a linear operator, can have infinitely many distinct zeros and each such $p_\infty$ will give rise to a different PDE for $h$. If we can compute $p_\infty$ then we can compute $p$ by computing $h$. In \cite{mandal2023learning} we provide a deep learning method to compute positive zeros of $\mathcal L$ and thus for our purposes it suffices to compute $h$. To compute $h$ with the Feynman-Kac formula we identify the following quantities, 

\begin{equation}
\begin{aligned}
    &\bar\mu = (\sigma^2\nabla\log p_\infty - \mu)\\
    &\bar\sigma =\sigma\\
    & q = 0
\end{aligned}
\end{equation}

In this construction $q$ has automatically become lower-bounded. Moreover, the formula \eqref{eq:FK} simplifies to 
\begin{align}
    h(t, \mathbf x) = \mathbb E\left[\left.\frac{p_0(\bar X_t)}{p_\infty(\bar X_t)}\right\vert \bar X_0=\mathbf x\right]\label{eq:h-E}
\end{align}
where $\bar X_t$ satisfies, 
\begin{align}
    d\bar X_t = (\sigma^2\nabla\log p_\infty-\mu)\,dt + \sigma\,dW_t\label{eq:h-SDE}
\end{align}
Finally, $p(t,\mathbf x)$ can be computed with the formula, 
\begin{align}
    p(t, \mathbf x) = \mathbb E\left[\left.\frac{p_0(X_t)}{p_\infty(X_t)}\right\vert X_0=\mathbf x\right]p_\infty(\mathbf x)\label{eq:p-E}
\end{align}

Note that even if the computed $p_\infty$ is unnormalized since we both multiply and divide by $p_\infty$ in \eqref{eq:p-E}, the normalization constant cancels out. Also, since $\bar\mu$ depends on $p_\infty$ through $\nabla\log p_\infty$, the normalization constant for $p_\infty$ is again eliminated and $X_t$ appearing in \eqref{eq:p-E} is independent of the normalization constant. Therefore, we can recover the correct $p$ with \eqref{eq:p-E} for any nonzero zero of $\mathcal L$ which can be computed following the recipe laid out in \cite{mandal2023learning}. For the sake of convenience, we will refer to equations \eqref{eq:h} and \eqref{eq:h-SDE}  as the \textit{h-equation} and the \textit{h-SDE} respectively.

We have only dealt with two of the three requirements of theorem~\ref{thm:FK-Oksendal} so far by appealing to the more general version of the Feynamn-Kac formula \cite{xi2019jump} and transforming the original FPE. The remaining requirement that the initial condition vanishes at infinity translates to $\frac{p_0}{p_\infty}$ vanishing at infinity. Even though we do not impose this restriction in our examples, in section~\ref{ssec-limit} we see that this condition is immaterial to our situation since we restrict the trajectories of the h-SDE \eqref{eq:h-SDE} to a finite domain and the expectation in \eqref{eq:p-E} is therefore bounded.

\label{ssec-Feynman-Kac}

\subsection{The algorithm}
Following the discussion in section~\ref{ssec-Feynman-Kac} we propose the following hybrid algorithm~\ref{algo:hybrid} for solving \eqref{eq:FPE-0} in high dimensions which uses the powers of both deep learning and the Feynman-Kac formula. In order to generate trajectories of the h-SDE \eqref{eq:ito-time-homo} we will use some time-discretized scheme like Euler-Maruyama or Milstein \cite{kloeden1992stochastic}. Just as in the prequel \cite{mandal2023learning}, we will restrict ourselves to computing the solution on a finite domain $\mathcal D$ where almost all the probability mass lies throughout $t\in[0, T]$. Besides the space-time boundaries $\mathcal D, T$, algorithm~\ref{algo:hybrid} has two other hyperparameters $M,N$ describing the time-discretization for Euler-Maruyama and the sample-size for estimating the Feynman-Kac expectation respectively.

\begin{algorithm}[!ht]
\begin{enumerate}
    \item Choose $\mathcal D, T, M, N$ according to the PDE and available computational resources.
    \item Use a known analytic zero of $\mathcal L$ or compute $p_\infty$ according to algorithm~6.1 in \cite{mandal2023learning}.
    \item Pick $t\in[0, T]$ and $\mathbf x\in\mathcal D$. Generate $N$ trajectories of the SDE given by, 
    \begin{align*}
        dX_\tau = (\sigma^2\nabla \log p_\infty - \mu)\,d\tau + \sigma\,dW_\tau
    \end{align*}
    starting from $\mathbf x$ and running till time $t$ with Euler-Maruyama method. Suppose the trajectories are computed at times $0=\tau_0<\tau_1<\cdots<\tau_M=t$. Let $X^{(i)}_{j}$ denote the $j$-th point (in time) in the $i$-th trajectory.
    \item Approximate $p(\tau_j, \mathbf x)$ with the following quantity, \begin{align*}
        \frac{p_\infty(\mathbf x)}{N}\sum_{i=1}^N\frac{p_0(X^{(i)}_{j})}{p_\infty(X^{(i)}_{j})}
    \end{align*}
    \end{enumerate}
    \caption{The $p_\infty$+FK algorithm}\label{algo:hybrid}
\end{algorithm}

\label{ssec-algo}

\subsection{Strengths and limitations}

The main strength of algorithm~\ref{algo:hybrid} lies in the fact that we can focus on the solutions pointwise without having to deal with global structures, thus somewhat mitigating the curse of dimensionality. Also, generating trajectories for the Feynman-Kac expectation with Euler-Maruyama lends itself to extreme parallelization, a must-have ingredient for solving high dimensional problems. The nature of the algorithm makes it applicable to a fairly large set of problems. Moreover, to compute solutions for the same problem but different initial conditions one can just reuse the same trajectory data and once the trajectory data are generated, solutions can be found for all of the time-points used to generate the trajectory data.

The limitations of algorithm~\ref{algo:hybrid} closely related to the properties dynamical system we are dealing with. In order to explore these, first, we need to explore the behavior of the h-SDE \eqref{eq:h-SDE} for different dynamical systems.
\subsubsection{Behavior of h-SDE trajectories and domain contraction}\label{ssec-h-behavior}
The stationary Fokker-Planck equation $\mathcal Lp_\infty=0$ has analytical solution for gradient systems. In this case, the solution can be expressed as 
\begin{align}
    &p_\infty\propto \exp\left(-\frac{2V}{\sigma^2}\right)\\
    \implies&\sigma^2\nabla\log p_{\infty} = 2\mu\\
    \implies&\bar{\mu}=\mu
\end{align}
where $V$ is defined as in \eqref{eq:grad-mu}, for a derivation see for example, section~3.1 in \cite{mandal2023learning}. Consequently, the h-SDE for gradient systems can be written as 
\begin{align}
    d\bar X_t = \mu\, dt+\sigma\, dW_t\label{eq:h-SDE-grad}
\end{align}
So for the gradient case, the h-SDE is identical to the original SDE \eqref{eq:SDE-0} describing the underlying dynamical system. Since the dynamical system possesses an attractor, the trajectories of the h-SDE are attracted to this this attractor and we can safely avoid blow-up even if the drift $\mu$ is non-Lipschitz.

But the same can not be said for the non-gradient systems where the drift term $\bar{\mu}$ in the h-SDE might be dominated by $-\mu$ rather than being $\mu$ as in the original SDE, thus losing its attracting properties. In such cases the solutions of the h-SDE might experience finite or infinite time blow-up.
\begin{figure}[!ht]
    \centering
\includegraphics[scale=0.55]{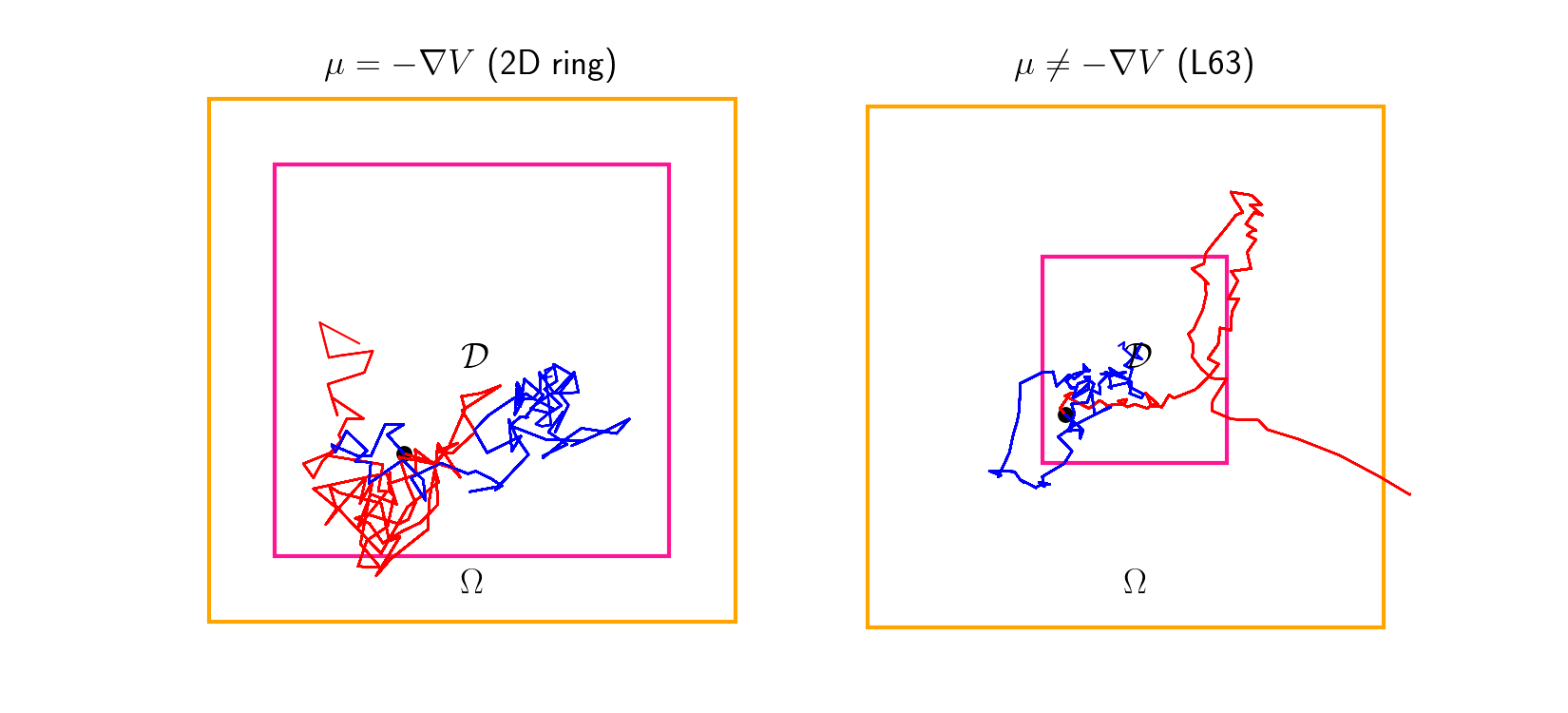}
    \caption{h-SDE trajectories for various systems. In both cases a pair of trajectories start from the same point (depicted as a black dot) in $\mathcal D$. While the trajectories for the gradient system might leave $\mathcal D$ (smaller rectangle), they do not leave $\Omega$ (larger rectangle). However, the same is not true for the non-gradient system.}
    \label{fig:h-SDE}
\end{figure}

Our knowledge of $\bar{\mu}$ is determined by our knowledge of $p_\infty$. If we do not have an analytic form for $p_\infty$ and have computed $p_\infty$ (up to the normalization constant) only up to the domain $\Omega$ then we can only be sure of our knowledge of $\bar{\mu}$ up to the domain $\Omega$. Consequently, if we want to compute the solution on the domain $\mathcal D$ then we must select $\mathcal D\subset\Omega$ in a way such that the trajectories of h-SDE that start inside $\mathcal D$, do not leave $\Omega$ till time $T$ with high probability. To quantify this notion, we can choose a tolerance $\varepsilon>0$ and select $\mathcal D, T$ such that,

\begin{align}
    \xi(T, \mathcal D, \Omega)\stackrel{\rm def}{=}\mathbb E_{\mathbf x\sim U(\mathcal D)}[P(\bar X_t\in\Omega\;\forall\;t\in[0, T]|\bar X_0=\mathbf x)] > 1-\varepsilon\label{eq:xi-defn}
\end{align}

\begin{figure}[!ht]
    \centering
\includegraphics[scale=0.5]{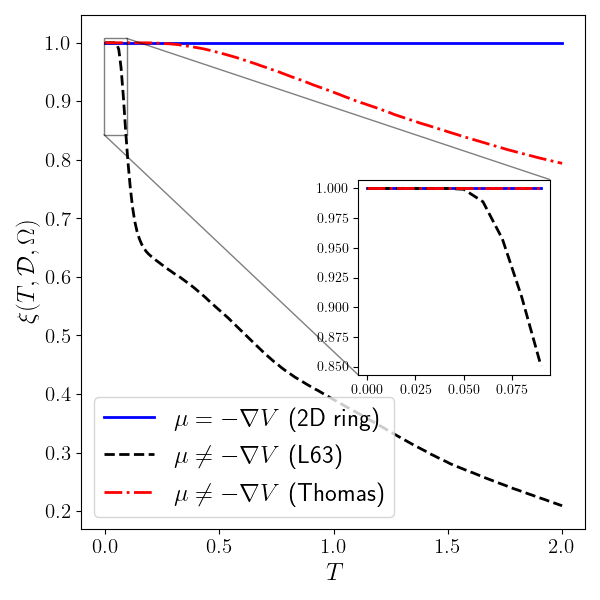}
    \caption{$\xi(T, \mathcal D, \Omega)$ as a function of $T$ for various systems.}
    \label{fig:xi}
\end{figure}

$\xi(T, \mathcal D, \Omega)$ denotes the average probability that a trajectory stays inside $\Omega$ till time $T$ given that it started inside $\mathcal D$, $U(\mathcal D)$ denotes the uniform distribution on $\mathcal D$ in \eqref{eq:xi-defn}. $\xi$ helps us specify the space-time boundaries for the effective employment of algorithm~\ref{algo:hybrid}. While for a gradient system, choosing $\Omega$ such that it contains the corresponding attractor, we can make sure that h-SDE trajectories originating from $\mathcal D$ do not leave $\Omega$ with high probability, the same can not be said for a non-gradient system. Figure~\ref{fig:h-SDE} shows the difference between h-SDE trajectories for typical gradient and non-gradient systems. Since for gradient systems the h-SDE trajectories do not leave $\Omega$ with high probability, $\xi(T, \mathcal D, \Omega)\approx 1$ for any $T$. In fact, empirically it might evaluate exactly to $1$. But for non-gradient system $\xi(T, \mathcal D, \Omega)$ is a decreasing function of $T$ as seen in figure~\ref{fig:xi} and in such cases we select the hyperparameter $T$ for algorithm~\ref{algo:hybrid} according to \eqref{eq:xi-defn} using a pre-chosen tolerance $\varepsilon$.

We are only able to solve \eqref{eq:FPE-0} on a subset $\mathcal D$ of where we solved its stationary counterpart, namely $\Omega$. We refer to this phenomenon as \textit{domain contraction}. Even though we have included gradient systems in the figures~\ref{fig:h-SDE}, \ref{fig:xi} for expository purposes, since we have perfect knowledge of $p_\infty$ (up to the normalization constant) over entire $\mathbb R^d$, domain contraction is a practical issue only for non-gradient systems and $T$ can be chosen arbitrarily large for gradient systems.

\subsubsection{Interpretation of finite time viability for non-gradient systems via Feynman-Kac on finite domains}\label{ssec-interpret-FK-finite}
Since we are dealing with finite domains, we can also look at the h-SDE through the lens of the Feynman-Kac formula on finite domains. In order to apply the finite domain version of the Feynman-Kac formula one requires perfect knowledge of the solution at the boundary at all times. In our case we need to know $h(t, \mathbf x)$ on $([0, T)\times{\partial\Omega}) \cup (\{0\}\times\bar{\Omega})$. For the ease of discussion let us define the following quantity,
\begin{align}
    h(t, \mathbf x) = \begin{cases}
        &\Psi(T, \mathbf x),\;\forall\;\mathbf x\in \bar{\Omega},\; t=0\\
        & \Psi(T-t, \mathbf x)\;\forall\;(t,\mathbf x)\in[0, T)\times\partial\Omega
    \end{cases}
\end{align}
Assuming we know $\Psi$, the Feynman-Kac formula for $h$ becomes,
\begin{align}
    h(t, \mathbf x) = \mathbb E\left[\left.\Psi(\tau\wedge T, \bar{X}_{\tau\wedge T})\right\vert \bar X_0=\mathbf x\right]\label{eq:h-E-finite}
\end{align}
where $\tau$ denotes the first exit-time of $\bar{X}$ for $\Omega$ or,
\begin{align}
    \tau(\mathbf x) \stackrel{\rm def}{=} \inf\{s>0: \bar X_s\not\in\Omega\} 
\end{align}
For derivations of such formulae the interested reader can refer to theorem~4.4.5 in \cite{gobet2016monte} or theorem~4.2 in chapter~7 of \cite{yong1999stochastic} and its preceding section. Since in reality we only know $h(t, \mathbf x)$ at $t=0$ and therefore have no knowledge of $\Psi(T-t, \mathbf x)$ at any point other than $t=0$, we can hope to effectively use \eqref{eq:h-E-finite} only if the trajectories of h-SDE do not exit $\Omega$ till time $T$ with high probability. In such a scenario $\tau\wedge T$ becomes equal to $T$ and we are not forced to evaluate $\Psi$ at any point where we have no knowledge of $\Psi$. Using the finite domain version of Feynman-Kac therefore leads us to similar conclusions as in the previous section, namely, algorithm~\ref{algo:hybrid} can be used to solve gradient and non-gradient systems up to arbitrarily large and finite times respectively.

\subsubsection{Effect allowing h-SDE trajectories to leave \texorpdfstring{$\Omega$}{Lg} in case of perfect knowledge of \texorpdfstring{$p_\infty$ in $\Omega$}{Lg}}\label{sssec-error}
 Since $p_\infty$ is not computed by algorithm~\ref{algo:hybrid} but rather is an input to it, a scenario of interest is when we have perfect (rather than approximate) knowledge of $p_\infty$ on $\Omega$. In such a scenario we would like the algorithm to produce reasonably good approximate solutions even when we are willing tolerate $\xi(T, \mathcal D, \Omega)\neq1$ or some h-SDE trajectories leaving $\Omega$ within our chosen $T$. To analyze this scenario let us define our knowledge of $p_\infty$ as,
 \begin{align}
     \hat p_\infty(\mathbf x) =\begin{cases}
         &p_\infty(\mathbf x), \qquad\mathbf x\in\Omega\\
         &p_\infty^\sharp(\mathbf x),\qquad\mathbf x\in\Omega^c
     \end{cases} 
 \end{align}
 where $p_\infty^\sharp\neq p_\infty$ represents imperfect knowledge of $p_\infty$ outside $\Omega$. Similarly we can define a modified h-SDE as,
 \begin{align}
     d\hat X_t = (\sigma^2\nabla\log \hat p_\infty -\mu)\,dt+\sigma\,dW_t
 \end{align} 
 and the corresponding solution generated by algorithm~\ref{algo:hybrid} as,
 \begin{align}
   &\hat h(t, \mathbf x) = \mathbb E\left[\left.\frac{p_0(\hat{X}_t)}{\hat{p}_\infty(\hat{X}_t)}\right\vert \hat{X}_0=\mathbf x\right]\\
   & \hat p(t ,\mathbf x)=\hat h(t, \mathbf x) \hat p_\infty(\mathbf x)
\end{align}
Now we are ready to analyze the pointwise error incurred when we allow h-SDE trajectories to escape $\Omega$ within time $T$ with probability less than $\varepsilon$ or $\xi(T, \mathcal D, \Omega)>1-\varepsilon$.
\begin{prop}Let $\bar E, \hat E$ be the events\footnote{Note that $\bar E, \hat E$ are events dependent on $t, \mathbf x$ but to avoid notational cluttering we do not make this dependence explicit.} that $\bar X, \hat X$ stay inside $\Omega$ till time $t$ after starting at $\mathbf x\in\mathcal D\subset\Omega$ respectively. Assume the following,
    \begin{enumerate}
    \item $\xi(T, \mathcal D, \Omega)>1-\varepsilon$
    
    \item $\exists$ a constant $K>0$ such that,
    \begin{align} 
     p_\infty(\mathbf x)\,\mathbb E[ h_0(\bar X_t)|\bar X_0=\mathbf x, \bar E^c],\;\hat p_\infty(\mathbf x)\,\mathbb E[\hat h_0(\hat X_t)|\hat X_0=\mathbf x, \hat E^c] < K \;\forall\;\mathbf (t,\mathbf x)\in[0, T]\times\mathcal D
    \end{align}
    where,\begin{align}
    &h_0 \stackrel{\rm def}{=} \frac{p_0}{p_\infty}\\
    &\hat h_0 \stackrel{\rm def}{=} \frac{p_0}{\hat p_\infty}
\end{align}
\end{enumerate}
Then, 
\begin{align}
    \mathbb E_{\mathbf x\sim U(\mathcal D)}[|\hat p(t, \mathbf x)-p(t, \mathbf x)|] <\, 2K\varepsilon\;\forall\;t\in[0, T]
    \label{eq:pt-error}
\end{align}\label{prop:escape}
\end{prop}
\begin{proof}
Let $\mathbf x\in\mathcal D$ and $\Delta_h = h_0-\hat h_0$.
    \begin{align}
    |\hat p(t, \mathbf x)-p(t, \mathbf x)| = |\Delta_1+\Delta_2|    
    \end{align}
where
\begin{align}
    \Delta_1 =  p_\infty(\mathbf x)\,\mathbb E\left[\left.\frac{p_0(\bar{X}_t)}{{p}_\infty(\bar{X}_t)}\right\vert \bar{X}_0=\mathbf x\right]-\hat p_\infty(\mathbf x)\,\mathbb E\left[\left.\frac{p_0(\bar{X}_t)}{\hat{p}_\infty(\bar{X}_t)}\right\vert \bar{X}_0=\mathbf x\right]
\end{align}
and,
\begin{align}
    \Delta_2 = \hat p_\infty(\mathbf x)\,\mathbb E\left[\left.\frac{p_0(\bar{X}_t)}{\hat{p}_\infty(\bar{X}_t)}\right\vert \bar{X}_0=\mathbf x\right]-\hat p_\infty(\mathbf x)\,\mathbb E\left[\left.\frac{p_0(\hat{X}_t)}{\hat{p}_\infty(\hat{X}_t)}\right\vert \hat{X}_0=\mathbf x\right]
\end{align}
Now note that,
\begin{equation}
\begin{aligned}
    \Delta_1 = &p_\infty(\mathbf x)\,\mathbb E[h_0(\bar X_t)|\bar X_0=\mathbf x, \bar E]\,P(\bar E)+p_\infty(\mathbf x)\,\mathbb E[h_0(\bar X_t)|\bar X_0=\mathbf x, \bar E^c]\,(1-P(\bar E))\\
    -&\hat p_\infty(\mathbf x)\,\mathbb E[\hat h_0(\bar X_t)|\bar X_0=\mathbf x, \bar E]\,P(\bar E)-\hat p_\infty(\mathbf x)\,\mathbb E[\hat h_0(\bar X_t)|\bar X_0=\mathbf x, \bar E^c]\,(1-P(\bar E))
\end{aligned}
\end{equation}
Recalling that we have perfect knowledge of $p_\infty$ inside $\Omega$ we get,
\begin{align}
    \Delta_1=&p_\infty(\mathbf x)\left[\mathbb E[\Delta_h(\bar X_t)|\bar X_0=\mathbf x, \bar E]\,P(\bar E)+\mathbb E[\Delta_h(\bar X_t)|\bar X_0=\mathbf x, \bar E^c]\,(1-P(\bar E))\right]\\
    =&p_\infty(\mathbf x)\,\mathbb E[\Delta_h(\bar X_t)|\bar X_0=\mathbf x, \bar E^c]\,(1-P(\bar E))
\end{align}
where we arrive at the last equality by noticing that $\Delta_h=0$ in the event of $\bar E$.
And,
\begin{equation}
\begin{aligned}
    \Delta_2 = &\hat p_\infty(\mathbf x)\,[\mathbb E[\hat h_0(\bar X_t)|\bar X_0=\mathbf x, \bar E]\,P(\bar E)+\mathbb E[\hat h_0(\bar X_t)|\bar X_0=\mathbf x, \bar E^c]\,(1-P(\bar E))]\\
    -&\hat p_\infty(\mathbf x)\,[\mathbb E[\hat h_0(\hat X_t)|\hat X_0=\mathbf x, \hat E]\,P(\hat E)+\mathbb E[\hat h_0(\hat X_t)|\hat X_0=\mathbf x, \hat E^c]\,(1-P(\hat E))]
\end{aligned}
\end{equation}
  Since $\bar X, \hat X$ follow identical dynamics inside $\Omega$, $P(\bar E)=P(\hat E)$ and we have,
\begin{align}
    \Delta_2 = &\hat p_\infty(\mathbf x)\,[\mathbb E[\hat h_0(\bar X_t)|\bar X_0=\mathbf x, \bar E^c]-\mathbb E[\hat h_0(\hat X_t)|\hat X_0=\mathbf x, \hat E^c]]\,(1-P(\hat E))
\end{align}
\end{proof}
Again noting $\hat p_\infty(\mathbf x)=p_\infty(\mathbf x)$ and $P(\bar E)=P(\hat E)$ we get,
\begin{align}
    \Delta_1+\Delta_2 =\, & [p_\infty(\mathbf x)\,\mathbb E[ h_0(\bar X_t)|\bar X_0=\mathbf x, \bar E^c]-\hat p_\infty(\mathbf x)\,\mathbb E[\hat h_0(\hat X_t)|\hat X_0=\mathbf x, \hat E^c]]\,(1-P(\hat E))\\
    \implies|\Delta_1+\Delta_2|<\,&2K\,(1-P(\hat E))
\end{align}
Since $\xi$ is a monotone decreasing function of its first argument, the first assumption is equivalent to,
\begin{align}
    \xi(t, \mathcal D, \Omega) > 1-\varepsilon\;\forall\; t\in[0, T]
\end{align}
Therefore, according to the definition of $\xi$,
\begin{align}
    \mathbb E_{\mathbf x\sim U(\mathcal D)}[P(\bar E)] > 1-\varepsilon \;\forall\;t\in[0, T]
\end{align}
which implies,
\begin{align}
    &\mathbb E_{\mathbf x\sim U(\mathcal D)}[1-P(\hat E)] <\varepsilon \;\forall\;t\in[0, T]\\
    \implies&\mathbb E_{\mathbf x\sim U(\mathcal D)}[|\Delta_1+\Delta_2|]<2K\varepsilon
\end{align}
which completes our proof.

The second assumption in proposition~\ref{prop:escape} makes sure if we compute the probability densities using only the trajectories that exit $\Omega$ before time $t$, we get bounded quantities independent of $\varepsilon$. Proposition~\ref{prop:escape} tells us\textit{, under suitable circumstances, the average error that is caused by allowing some h-SDE trajectories to leave $\Omega$ is proportional to the fraction of trajectories that leave $\Omega$.} If we modify the first assumption to require that the inequality holds pointwise for every $\mathbf x\in\mathcal D$ instead, we can bound the supremum norm of the error over $\mathcal D$ instead of the average error.
\label{ssec-limit}

\section{Results}
\label{sec-results}
\begin{figure}[!ht]
    \centering
\includegraphics[scale=0.41]{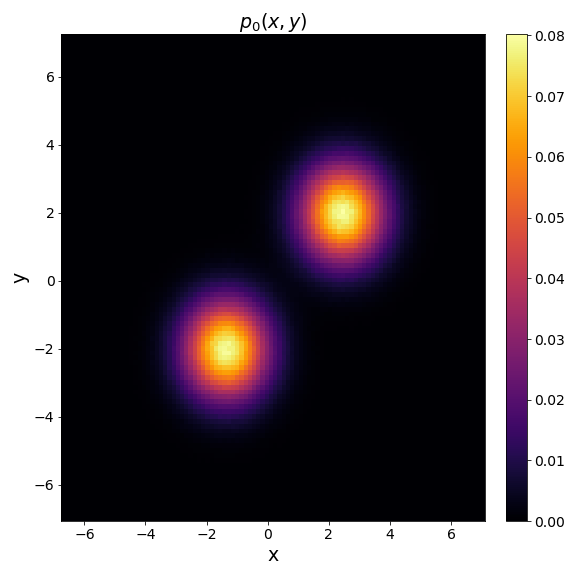}
    \caption{Initial condition for the noisy Lorenz-63 system.}
    \label{fig:L63-0}
\end{figure}

In this section we describe the results in a manner that parallels the examples in section~\ref{sec-example} with additional problem-specific details. We refer to the solutions obtained via algorithm~\ref{algo:hybrid} as the \textit{learned} solutions in the figures below. We use the Monte-Carlo method as detailed in algorithm~\ref{algo:MC} to produce reference solutions in 2 and 3 dimensions. But in the absence of analytical solutions, we refrain from calculating the difference between the reference solution and the solution produced by algorithm~\ref{algo:hybrid} since Monte-Carlo solutions can be significantly erroneous, sometimes by more than an order of magnitude compared to its counterpart, even for simple problems, as evidenced by figure~7.2 in the prequel \cite{mandal2023learning}. When dealing with high-dimensional PDEs, one widely accepted paradigm is to construct solutions that are \textit{statistically} accurate or have the correct coarse-grained structures, especially for particle-based methods in geophysics \cite{bosler2013particle}, \cite{chen2018efficient} and financial modelling \cite{cui2015particle}. We present our results in the same spirit. For each system we use a bimodal initial condition as described in section~\ref{sec-example}. Figure~\ref{fig:L63-0} shows the initial condition for the noisy Lorenz-63 system defined in \eqref{eq:p0-L63}. For ease of visualization we have integrated out the last coordinate. The symmetry of the initial condition implies the other 2D marginals are identical. The other initial conditions described in section~\ref{sec-example} are either identical or qualitatively similar. In order to produce 2D marginal densities we use Gauss-Legendre quadrature the details of which can be found in appendix~9.3 of \cite{mandal2023learning}.

\subsection{10D ring system} 
\begin{figure}[!ht]
    \centering
\includegraphics[scale=0.32]{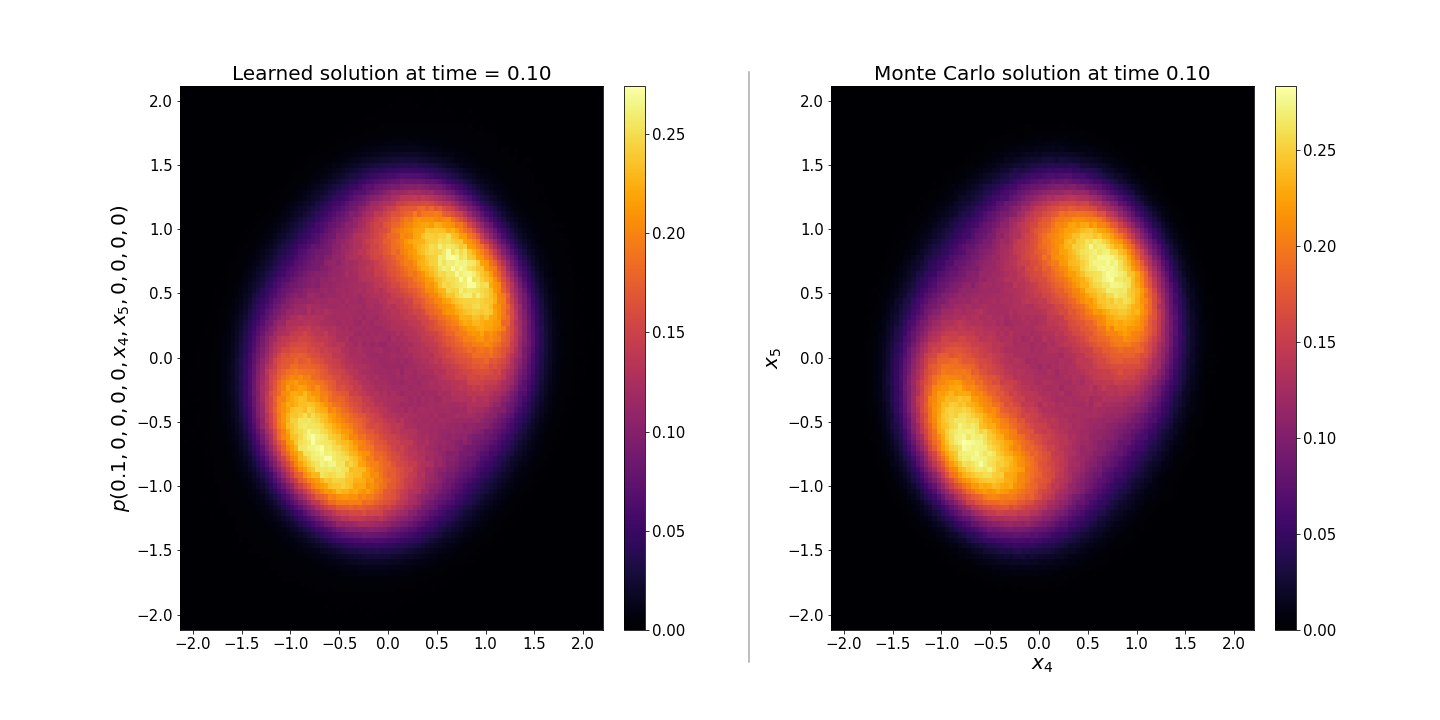}
    \caption{Solutions for the 10D time-dependent system at time $t=0.1$. The learned solution has been normalized such that $\int_\mathbb{R}\int_\mathbb{R} p(0.1, 0, 0, 0, 0, x_4, x_5, 0, 0, 0, 0)\, dx_4\,dx_5=1$. The right panel depicts the Monte-Carlo solution for the $2$D Fokker-Planck equation corresponding to the variables $x_4,x_5$. The learned and Monte-Carlo solutions were computed using $10^5$ (pointwise) and $10^7$ trajectories respectively.}
    \label{fig:10D-time}
\end{figure}

Figure~\ref{fig:10D-time} shows the computed solution for the 2$n$D ring system presented in section~\ref{ssec-2nD} for $n=5$, at time $t=0.1$. More specifically, the left panel of figure~\ref{fig:10D-time} depicts the probability density when all but the variables $x_4, x_5$ are set to $0$. For better visualization and comparison with a reference solution, the \textit{learned} solution has been normalized in a way such that,
\begin{align*}
   \int_\mathbb{R}\int_\mathbb{R} p(0.1, 0, 0, 0, 0, x_4, x_5, 0, 0, 0, 0)\, dx_4\,dx_5=1 
\end{align*}
Here we have used $M=10$, $N=10^5$ and $\mathcal D=[-2,2]^{10}$ for algorithm~\ref{algo:hybrid}.

There are no analytical solutions for this problem and classical methods are unsuitable due to the curse of dimensionality. Therefore, in order to compute a reference solution we use the careful design of the problem itself. The system described by \eqref{eq:potential-2nD} and the initial condition \eqref{eq:p0-2nD} is essentially $n$ identical, decoupled $2$D systems and we can approximately solve this $2$D system with Monte-Carlo method which is depicted in the right panel of figure~\ref{fig:10D-time}. We use $10^7$ trajectories to generate the Monte-Carlo solution. Note that, to compute the \textit{learned} solution in figure~\ref{fig:10D-time} we use a neural network approximation of $p_\infty$ obtained via the method described in \cite{mandal2023learning} rather than using the analytical version of $p_\infty$. This demonstrates the viability of algorithm~\ref{algo:hybrid} in conjunction with algorithm~6.1 of \cite{mandal2023learning} for high dimensional problems.

\subsection{Noisy Lorenz-63 system}\label{ssec-res-L63}
\begin{figure}[!ht]
    \centering
\includegraphics[scale=0.21]{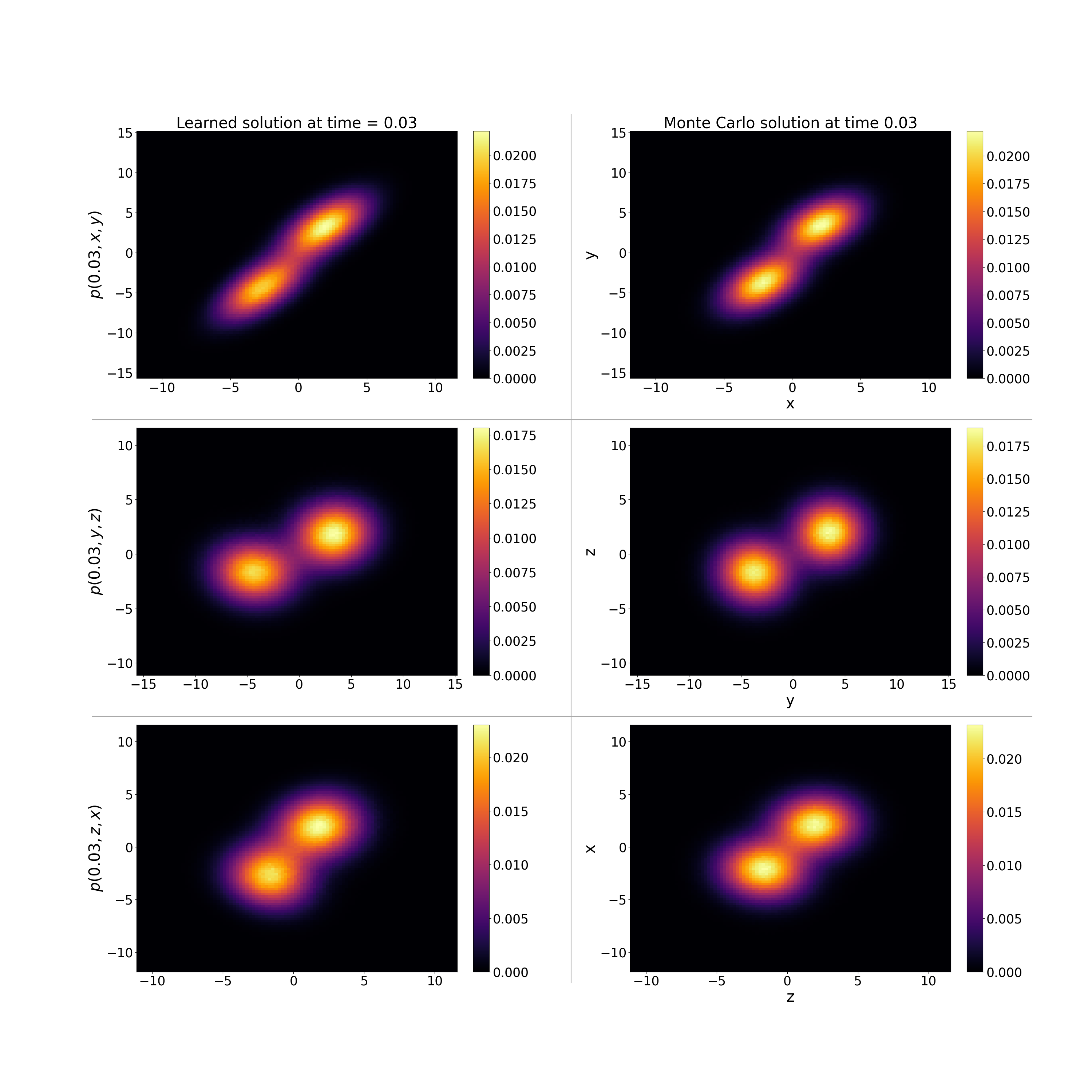}
    \caption{Solutions for the noisy Lorenz 63 system at time t=0.03. The learned and Monte-Carlo solutions were computed using $200$ (pointwise) and $10^7$ trajectories respectively.}
    \label{fig:L63-time}
\end{figure}

Figure~\ref{fig:L63-time} shows the solutions for the noisy Lorenz system defined by \eqref{eq:mu-L63} at time $t=0.03$ with $\mathcal D=[-10, 10]\times[-15, 15]\times[7, 28], \Omega=[-30, 30]\times[-40, 40]\times[0, 70]$. As seen in figure~\ref{fig:xi}, $\xi(t,  \mathcal D, \Omega)\approx 1$. For easier visualization we present the 2D marginals $p(t, x, y), p(t, y, z), p(t, z, x)$. To compute $p_\infty$ in a functional form for this problem we use algorithm~6.1 in \cite{mandal2023learning}. We use $M=3$ and $N=200$ for the learned solution and $10^7$ trajectories to generate the corresponding Monte-Carlo solution. This shows that we can produce solutions with algorithm~\ref{algo:hybrid} that are comparable to Monte-Carlo method using 
several orders of magnitude fewer trajectories for each point. Both methods use $0.01$ as the step length for Euler-Maruyama.

\subsection{Noisy Thomas system}
\begin{figure}[!ht]
    \centering
\includegraphics[scale=0.32]{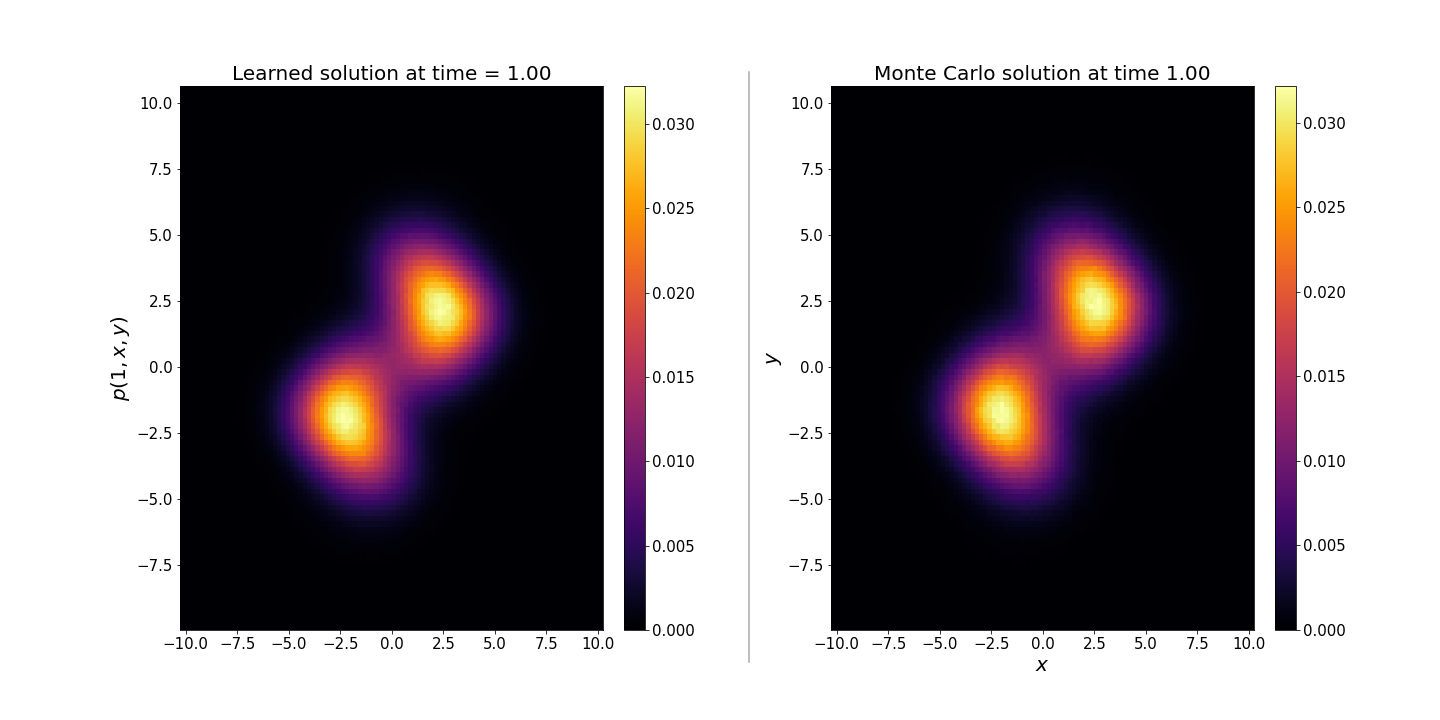}
    \caption{Solutions for the noisy Thomas system at time t=1. The learned and Monte-Carlo solutions were computed using $50$ (pointwise) and $10^7$ trajectories respectively.}
    \label{fig:Thomas-time}
\end{figure}

Figure~\ref{fig:Thomas-time} shows the solutions for the noisy Thomas system defined by \eqref{eq:mu-Thomas} at time $t=1.0$ with $\mathcal D=[-8, 8]^3$ and $\Omega=[-10, 10]^3$. We only present $p(t, x, y)$, noting that the symmetry of the problem renders demonstrations of the other 2D marginals redundant. We use algorithm~6.1 in \cite{mandal2023learning} for computing $p_\infty$. We use $M=10, N=50$ for the learned solution and the Monte-Carlo counterpart is computed using $10^7$ trajectories which reinstates our intuition that computing the Feynman-Kac expectation requires far fewer trajectories compared to Monte-Carlo for a similar level of accuracy.  Both methods use $0.1$ as the step length for Euler-Maruyama. In figure~\ref{fig:xi} we see that $\xi(t, \mathcal D, \Omega)\approx0.916$. So even after letting nearly $8.4\%$ of the h-SDE trajectories escape $\Omega$, we achieve a reasonable approximation.

\subsection{One step filter}
\begin{figure}[!ht]
    \centering
\includegraphics[scale=0.21]{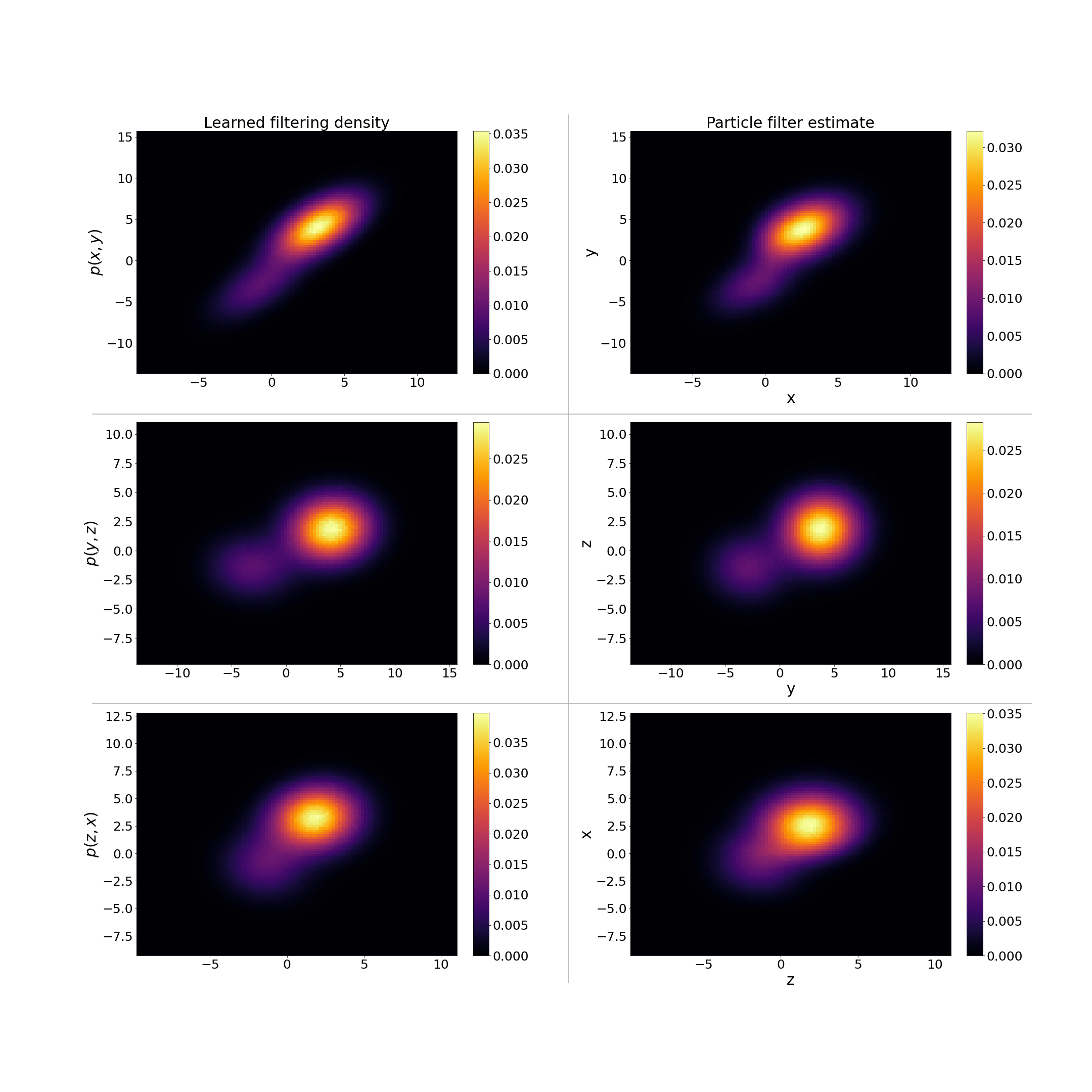}
    \caption{One step filtering density for the noisy Lorenz 63 system. The learned and particle filter solutions were computed using $200$ (pointwise) and $10^7$ trajectories respectively.}
    \label{fig:L63-filter}
\end{figure}

Suppose in the filtering problem described in section~\ref{ssec-1-filter-prob} we only observe $x,z$ coordinates of the noisy Lorenz 63 system and our observation noise standard deviation is $\sigma_o=5$. The one-step filtering density can be written as 
\begin{align}
    p(x_1|y_1) \propto p(x_1|x_0)p(y_1|x_1)
\end{align}
For a derivation see chapter 6 of \cite{sarkka2023bayesian} or \cite{doucet2009tutorial}. $p(x_1| x_0)$ can be thought of as the solution to the corresponding Fokker-Planck equation at time $t=g=0.03$ (the observation gap) with the initial condition being equal to the density of $x_0$,. We can calculate the likelihood $p(y_1|x_1)$ using $\sigma_o$ which lets us estimate the 2D marginals of the one-step filtering density for this problem. We calculate the solution to the Fokker-Planck equation with algorithm~\ref{algo:hybrid} and Monte-Carlo. The final one-step filtering densities are shown in figure~\ref{fig:L63-filter}. Computing the filtering density with Monte-Carlo in this way is akin to using the bootstrap particle filter, a popular nonlinear filtering algorithm, see chapter 11 of \cite{sarkka2023bayesian} or \cite{doucet2009tutorial} for more discussion on particle filters. Therefore, we refer to the Monte-Carlo estimate for the filtering density as the particle filter estimate in figure~\ref{fig:L63-filter}. We use the same $M, N$ and the same number of trajectories for Monte-Carlo as we did in section~\ref{ssec-res-L63}. It is interesting to note that the solution in figure~\ref{fig:L63-time} is bimodal whereas in the filtering density in figure~\ref{fig:L63-filter} one of the mode collapses after we make an observation.

\section{Conclusions and future work}
\label{sec-conclusions}
Physics-informed neural networks are unsuitable for a large class of partial differential equations of parabolic type. But algorithm~\ref{algo:hybrid} in conjunction with algorithm~6.1 in \cite{mandal2023learning} gives us a viable method for computing solutions to Fokker-Planck equations in high dimensions. The ability to focus on the solution pointwise lets us avoid the curse of dimensionality in a practical sense. To compute the Feynman-Kac expectaion we can compute the Euler-Maruyama trajectories in a highly parallel way and the same trajectories can be used to calculate solutions for different initial conditions. We are able to compute solutions for gradient systems up to arbitrarily long times. But for non-gradient system we encounter domain contraction, a phenomenon where we are able to calculate the solution on a subset of where know the stationary solution. Due to blow-up of the h-SDE for non-gradient systems we can compute solutions only up to a finite time. In case we have perfect knowledge of $p_\infty$, under ideal conditions the error that is introduced by letting some trajectories of h-SDE escape the region where we know the stationary solution is proportional to the fraction of trajectories that escape. Even after allowing a significant portion of trajectories to escape, we can achieve decent approximations. Moreover, algorithm~\ref{algo:hybrid} can achieve solutions that are comparable to Monte-Carlo solutions using orders of magnitude fewer trajectories for each point. The connection between the Fokker-Planck equation and the stochastic filtering problem suggests that modifications/extensions of algorithm~\ref{algo:hybrid} can be useful for devising new filtering algorithms.

\section*{Acknowledgements}
This work was supported by the Department of Atomic Energy, Government of India, under project no. RTI4001.

\section{Appendix}
\label{sec-appendix}
\subsection{Existence and uniqueness of solutions to example problems}\label{ssec-unique} Since appropriately scaling $t,\mu$ we can change \eqref{eq:FPE-1} to have unit diffusion, according to section 1 of \cite{bogachev2021uniqueness} and theorem 9.4.6 and example 9.4.7 of \cite{bogachev2022fokker} it suffices to prove that $\exists\, C>0$ such that
\begin{align}
    \mu\cdot \mathbf x \le C+C\|\mathbf x\|_2^2\,,\qquad\forall\;\mathbf x\in\mathbb R^d
\end{align}
in order to confirm the existence of a unique solution to \eqref{eq:FPE-0}. 

It is easy to check that the above condition is satisfied for the drift terms for each of the examples discussed in section~\ref{sec-example}. For the 2D ring system in~\eqref{eq:ring2D}, we note that $\mu\cdot\mathbf x= -4r^2(r^2-1)\le 1$, using the fact that the function $f(x)=x(1-x)$ achieves global maximum at $x=\frac{1}{2}$, so we can choose $C = 1$. For the noisy L63 system defined by \eqref{eq:mu-L63} we have $\mu\cdot\mathbf x =-\alpha x^2-y^2-\beta z^2+(\alpha+\rho) xy\le (\alpha+\rho)r^2$ so it suffices to set $C=(\alpha+\rho)$. Lastly, for the noisy Thomas system defined by \eqref{eq:mu-Thomas} we have $\mu\cdot\mathbf x =-br^2 + x\sin y + y\sin z + z\sin x \le 3r\le 3+3r^2$, therefore $C=3$ is a suitable choice.

\subsection{Monte Carlo algorithm}\label{ssec-algo-MC}
The relationship between \eqref{eq:SDE-0} and \eqref{eq:FPE-0}  \cite{risken1996fokker}, \cite{bogachev2022fokker} gives us the following way of estimating solutions of Fokker-Planck equations. We can evolve multiple particles according to \eqref{eq:SDE-0} up to time $t$ using Euler-Maruyama method \cite{kloeden1992stochastic}, subdivide the domain $\mathcal D$ into $d$-dimensional boxes and count the how many particles lie inside each box to compute the probability density at the centers of these boxes. Here $\mathcal N$ denotes the multivariate normal distribution.
\begin{algorithm}[!htp]
Sample $\{ X_0^{(i)}\}_{i=1}^N\sim p_0$.\\
Set the time-step $h=\frac{t}{M}$.\\
\For {$j=1,2\cdots, M$}{
 Sample $w^i_k\sim\mathcal N(\mathbf 0_d, h I_d)\;\;\forall\;i$\\
 $ X_k^{(i)}\leftarrow  X_{k-1}^{(i)} + \mu\left(X_{k-1}^{(i)}\right)h + \sigma w_k^i\;\;\forall\;i$\\
}
Subdivide the domain of interest $\mathcal D$ into $d$-dimensional boxes.\\ Count the number of $X^{(i)}_{S}$ that are in a box to estimate the probability density at the center of the box.
\caption{Monte Carlo algorithm}\label{algo:MC}
\end{algorithm}

\bibliographystyle{siamplain}
\bibliography{ref}

\end{document}